\documentclass[]{article}
\usepackage[margin=1in]{geometry}
\usepackage{amsfonts,amsmath}
\usepackage{tabularx, svg}
\usepackage{hyperref}
\usepackage{graphicx}
\usepackage{bigfoot}
\usepackage[norelsize, linesnumbered, ruled, lined, boxed, commentsnumbered]{algorithm2e}
\title{Simple exponential acceleration of the power iteration algorithm}
\author{Congzhou M Sha$^{a,b}$  (\url{cms6712@psu.edu}), Nikolay V Dokholyan$^{a,b,c,d,*}$ (\url{dokh@psu.edu})}
\date{%
	$^a$Department of Pharmacology, Penn State College of Medicine\\
	700 HMC Cres Rd, Hershey, PA, United States\\
	$^b$Department of Engineering Science and Mechanics, Penn State University\\
	212 Earth-Engineering Sciences Bldg., State College, PA, United States\\
	$^c$Department of Biochemistry and Molecular Biology, Penn State College of Medicine\\
	700 HMC Cres Rd, Hershey, PA, United States\\
	$^d$Departments of Chemistry and Biomedical Engineering, Penn State University\\
	77 Pollock Rd, State College, PA, United States\\
	*Corresponding author\\[2ex]
	\today
}

\begin{document}

\maketitle

\begin{abstract}
Many real-world problems rely on finding eigenvalues and eigenvectors of a matrix. The power iteration algorithm is a simple method for determining the largest eigenvalue and associated eigenvector of a general matrix. This algorithm relies on the idea that repeated multiplication of a randomly chosen vector $x$ by the matrix $A$ gradually amplifies the component of the vector along the eigenvector of the largest eigenvalue of $A$ while suppressing all other components. Unfortunately, the power iteration algorithm may demonstrate slow convergence. In this report, we demonstrate an exponential speed up in convergence of the power iteration algorithm with only a polynomial increase in computation by taking advantage of the commutativity of matrix multiplication.
\end{abstract}
\pagebreak
\section{Introduction}
\par Eigendecomposition of matrices is useful in continuum mechanics \cite{landaulifshitz}, quantum mechanics \cite{sakurai}, machine learning \cite{pca}, graph theory \cite{spectralgraph}, and many other areas of science and mathematics.
\par The power iteration algorithm is a well-known method for determining the largest eigenvalue and its associated eigenvector of a general matrix \cite{original}. The following argument is taken from \cite{bindellecture}. Given a diagonalizable matrix $A\in \mathbb{C}^{n\times n}$, we may write
\begin{equation} \label{eq:1}
	A=U\Lambda U^{-1}
\end{equation}
\begin{equation} \label{eq:2}
	A^k=(U\Lambda U^{-1})(U\Lambda U^{-1})\cdots (U\Lambda U^{-1})=U\Lambda^kU^{-1}
\end{equation}
where the columns of $U$ are the eigenvectors $v_i$, and $\Lambda$ is diagonal with the eigenvalues $\lambda_i$ of $A$ along the diagonal. Multiplying Equation \ref{eq:2} on the right by $U$,
\begin{equation}
	A^kU=U\Lambda^k
\end{equation}
Randomly choosing a vector $x=Uy\in \mathbb{C}^n$, we may write
\begin{equation}
	A^kx=U\Lambda^ky=\lambda_1^k\left(\sum_{i=1}^nv_i\left(\frac{\lambda_i}{\lambda_1}\right)^ky_i\right)
\end{equation}
Here,  we assume that $\lambda_1$ is the largest eigenvalue of $A$, and $v_1$ is its associated eigenvector. With each factor of $A$, the components along the other eigenvectors are suppressed by a factor of $|\frac{\lambda_i}{\lambda_1}|<1$, and thus we expect $A^kx$ to become increasingly parallel to $v_1$.
\par The power iteration algorithm is usually expressed as below, where $x_0$ is drawn randomly and $f$ indicates the final value of $k$:\\
\begin{algorithm}[H]
\caption{Power Iteration}\label{alg:pi}
	\SetAlgoLined
	\LinesNumbered
	\SetKwInOut{Input}{Input}
	\Input{$A\in\mathbb{R}^{n\times n}$, tolerance $tol$}
	$x_0\in\mathbb{R}^n$\\
	$k\leftarrow 0$\\
	\SetKwRepeat{Do}{do}{while}
	\Do{$|x_{k+1}-x_{k}|>tol$}{ 
		$x_{k+1}\leftarrow \left(Ax_{k}\right)/\left(|Ax_{kz}|\right)$\\
		$k\leftarrow k+1$}
	\Return $(x_{f}^\dagger Ax_{f})/(x_{f}^\dagger x_{f})$ and $x_f$\\
	\SetKwInOut{Output}{Output}
	\Output{Rayleigh quotient estimate of the largest eigenvalue $\lambda$, estimated eigenvector $x$}
\end{algorithm}
However, there is a well-known algorithm for the computation of high powers of a matrix: exponentiation through squaring. Without loss of generality, we assume that $k=2^j, j\in \mathbb{N}$. $f$ again indicates the final value of $i$.\\
\begin{algorithm}[H]
	\caption{Matrix exponentiation through squaring}\label{alg:exp}
	\SetAlgoLined
	\LinesNumbered
	\SetKwInOut{Input}{Input}
	\Input{$A\in\mathbb{R}^{n\times n}$, desired power $k=2^j$}
	$A_0\leftarrow A$\\
	$i\leftarrow 0$\\
	\SetKwRepeat{Do}{do}{while}
	\Do{$i\leq j$}{ 
		$A_{i+1}\leftarrow A_i\cdot A_i$\\$i\leftarrow i+1$}
	\Return $A_i$
	\SetKwInOut{Output}{Output}
	\Output{$A_f=A^k$}
\end{algorithm}
\section{Material and Methods}
\subsection{Proposed algorithms}
We propose an algorithm which combines Algorithms \ref{alg:pi} and \ref{alg:exp} to achieve exponential convergence of the largest eigenvalue:\\
\begin{algorithm}[H]
	\caption{Power iteration with exponentiation}\label{alg:pie}
	\SetAlgoLined
	\LinesNumbered
	\SetKwInOut{Input}{Input}
	\Input{$A\in\mathbb{R}^{n\times n}$, tolerance $tol$}
	$x_0\in\mathbb{R}^n$\\
	$A_0\leftarrow A$\\
	$i\leftarrow 0$\\
	\SetKwRepeat{Do}{do}{while}
	\Do{$|A_{i+1}-A_{i}|>tol$}{ 
		$A_{i+1}\leftarrow A_i\cdot A_i$\\$A_{i+1}\leftarrow A_{i+1}/|A_{i+1}|$\\$i\leftarrow i+1$}
	$x_f\leftarrow A_fx_0$\\
	\SetKwInOut{Output}{Output}
	\Return $(x_{f}^\dagger A_0x_{f})/(x_{f}^\dagger x_{f})$ and $x_f$\\
	\Output{Rayleigh quotient estimate of the largest eigenvalue $\lambda$, estimated eigenvector $x$}
\end{algorithm}
\par To generalize Algorithms \ref{alg:pi} and \ref{alg:pie} to the problem of finding the largest $k$ eigenvectors in the case of a self-adjoint (Hermitian) matrix, we may perform a version of Gram-Schmidt orthonormalization carried out on the input matrix $A$, as shown in Algorithm \ref{alg:gs}.\\
\begin{algorithm}[H]
	\caption{Gram-Schmidt Orthonormalization}\label{alg:gs}
	\SetAlgoLined
	\LinesNumbered
	\SetKwInOut{Input}{Input}
	\Input{$A\in\mathbb{R}^{n\times n}$, desired number of largest eigenvalues $k$}
	$A_0\leftarrow A$\\
	$i\leftarrow 0$\\
	$B\leftarrow\{e_i\}$\\
	$\Lambda\leftarrow [\ ]$\\
	\SetKwRepeat{Do}{do}{while}
	\Do{$i\leq k$}{
		$\lambda, v\leftarrow \text{PowerIteration}(A_i)$\\
		Append $\lambda$ to $\Lambda$\\
		$e_i\leftarrow v/|v|$\\
		$j\leftarrow 0$\\
		$A_{i+1}\leftarrow A_{i}-\lambda e_ie_i^\dagger$
	}
	\Return $B$, $\Lambda$\\
	\SetKwInOut{Output}{Output}
	\Output{The largest $k$ eigenvectors and eigenvalues of $A$}
\end{algorithm}
To demonstrate Algorithm \ref{alg:gs}, assume we have computed the largest eigenvalue and eigenvector of $A$ in line 6 of Algorithm \ref{alg:gs} (by using either Algorithm \ref{alg:pi} or \ref{alg:pie}). After normalizing the eigenvector (line 8), Then we can remove that component from the matrix $A$ in line 10
\begin{equation} \label{eq:4}
	(A-\lambda_1v_1v_1^\dagger)v_1=\lambda_1v_1-\lambda_1v_1(v_1^\dagger v_1)=0
\end{equation}
Thus $v_1$ now has eigenvalue $0$ for the matrix $A-\lambda_1v_1v_1^T$. Because we have the additional assumption that $A$ is self-adjoint ($A=A^\dagger$), its eigendecomposition is guaranteed to be unitary, by the finite-dimensional spectral theorem \cite{st}. Therefore, the Gram-Schmidt process works since at each step, we do not change the span of the remaining eigenvectors.
\subsection{Implementation and computational resources}
\par We implemented these algorithms and visualized results in Python 3.8 \cite{python}, used the NumPy 1.19 package for linear algebra as reference \cite{numpy}, and created figures with Matplotlib 3.4 \cite{matplotlib}. The code is included in the supplemental data (\verb|*.py|). Matrices of various dimensions were generated randomly with entries drawn from the normal distribution (mean 0, standard deviation 1) and symmetrized, to guarantee real eigenvalues and eigenvectors in the case of real matrices. The largest eigenvalue of each matrix was computed using the NumPy linear algebra library (\verb|numpy.linalg.eigvals|). The matrix norms in Algorithms \ref{alg:pi} line 4 and Algorithm  \ref{alg:pie} line 6 are necessary to prevent numerical overflow, but are otherwise arbitrary, so we chose the efficient max norm (maximum matrix element).
\par We recorded the computation time and the number of iterations of Algorithms \ref{alg:pi} and \ref{alg:pie} required to reach a numerical tolerance of $10^{-10}$. Calculations were performed on 12 cores of an AMD Ryzen 9 3900X CPU, with 64 GB of RAM available, allowing for parallelized matrix-matrix and matrix-vector multiplication through the NumPy library.
\section{Results}
We observed significant speedups in using Algorithm \ref{alg:pie} as compared to Algorithm \ref{alg:pi} for both real and complex matrices, shown in Tables 1 and 2\footnote{Table 1 corresponds to the file \verb|test.py| and Table 2  to \verb*|testcomplex.py|.}. We also observed that it took far fewer iterations to reach numerical convergence using Algorithm \ref{alg:pie} compared to Algorithm \ref{alg:pi}. We observed several instances in which it took Algorithm \ref{alg:pi} over $10^6$ iterations to reach the desired tolerance, which would occasionally significantly skewed the timings ($n=100$, 300 matrices taking almost 60 s to complete). In contrast, Algorithm \ref{alg:pie} consistently reached convergence within 20 iterations, shown in Figure 1.
\begin{center}
\begin{table}[h]
	\centering
	\label{table:bench}
	\caption{Benchmarking results (real matrices)}
	\begin{tabular}{|l|l|l|l|l|}
		\hline&
		$n=100$ (300 matrices) &
		$n=1000$ (5 matrices) &
		$n=3000$ (1 matrix) &
		$n=5000$ (1 matrix) \\
		\hline
		Alg 1 &
		35 s (117 ms/matrix) &
		29 s (5.8 s/matrix) &
		59 s &
		110 s \\
		Alg 3 &
		0.54 s (1.8 ms/matrix) &
		0.77 s (153 ms/matrix) &
		2.8 s &
		9.7 s \\\hline
		speedup & 65$\times$ & 36$\times$ & 21$\times$ & 11$\times$\\
		\hline
	\end{tabular}
\end{table}
\end{center}
\begin{center}
	\begin{table}[h]
		\centering
		\label{table:bench2}
		\caption{Benchmarking results (complex matrices)}
		\begin{tabular}{|l|l|l|l|l|}
			\hline&
			$n=100$ (300 matrices) &
			$n=1000$ (5 matrices) &
			$n=3000$ (1 matrix) &
			$n=5000$ (1 matrix) \\
			\hline
			Alg 1 &
			67 s (223 ms/matrix) &
			69 s (14 s/matrix) &
			175 s &
			359 s \\
			Alg 3 &
			1.38 s (4.6 ms/matrix) &
			7.3 s (1.5 s/matrix) &
			25 s &
			95 s \\\hline
			speedup & 49$\times$ & 9.5$\times$ & 7$\times$ & 3.8$\times$\\
			\hline
		\end{tabular}
	\end{table}
\end{center}
Even though Algorithm \ref{alg:pie} was written in Python whereas NumPy code executes highly optimized LAPACK routines for eigendecomposition, we also observed that Algorithm \ref{alg:pie} (taken advantage of LAPACK matrix multiplication) is faster than NumPy's \verb|numpy.linalg.eigvals| when we desire only the largest eigenvalue (3$\times$ speed up on three hundred $100\times 100$ matrices). It may be beneficial to write a version of Algorithm \ref{alg:pie} in a lower level language like C or Fortran to perform a more direct comparison.
\section{Discussion}
Algorithm \ref{alg:pie} computes a high (normalized) power of $A$, which is used to project $x_0$ to $x_f$, which as mentioned in the introduction should be close to $v_1$. Crucially, although naive multiplication of $n\times n$ matrices in line 5 of Algorithm \ref{alg:pie} takes $O(n^3)$ time, the powers $\left(\frac{\lambda_i}{\lambda_1}\right)^k$ take $O(\lg k)$ steps to compute and thus there is an exponential increase in accuracy with each matrix multiplication performed. Overall, Algorithm \ref{alg:pie} takes $O(n^3\lg k)$ time, using naive matrix multiplication. In contrast, Algorithm \ref{alg:pi} takes $O(k)$ steps to compute $\left(\frac{\lambda_i}{\lambda_1}\right)^k$, and thus $O(n^2 k)$ time. Thus there is a tradeoff between the cost of matrix multiplication and desired accuracy, in which Algorithm \ref{alg:pie} appears to be more performant in moderately-sized practical computations ($n=100$).
\par Just as in the traditional power iteration algorithm, Algorithm \ref{alg:pie} relies on the choice of an appropriate $x_0$. Random choice of $x_0$ should prevent the degenerate case of $x_0$ orthogonal to $v_1$ except in rare cases. Algorithm \ref{alg:pie} may also be more efficient with sparse matrices, as the cost of $O(n^3)$ matrix-matrix multiplication may be greatly reduced, though further benchmarking is required to compare with sparse applications of Algorithm \ref{alg:pie}. For other applications, the throughput of GPU-accelerated matrix-matrix multiplication may further equalize the performance of matrix-matrix multiplication in Algorithm \ref{alg:pie} compared to the matrix-vector multiplication in Algorithm \ref{alg:pi}.
\par The major drawback of Algorithm \ref{alg:pie} is the same as in traditional power iteration: only the largest eigenvalue and associated eigenvector are computed. Other iterative methods such as the Lanczos algorithm \cite{lanczosorig, lanczos} and inverse iteration \cite{inverseit} may compute other eigenvalues. It may also be possible to incorporate Algorithm \ref{alg:exp} into these iterative methods to speed up convergence. We also provided a method for computing the $k$ largest eigenvalues and associated eigenvectors in Algorithm \ref{alg:gs} when applied to self-adjoint matrices, which we observed to work for small values of $k$ (see provided code). Algorithm \ref{alg:gs}\footnote{Toy code is shown in \verb|eigenbasis.py|.} may be able to compete with other iterative methods for small $k$ and large $n$, though benchmarking is needed.
\par Even with its drawbacks, Algorithm \ref{alg:pie} may have practical uses. The largest eigenvalue and associated eigenvector can be highly informative, for example in the graph similarity score (Equation 5) of Zager and Verghese \cite{graphsim}. There may be other specific use cases where Algorithm \ref{alg:pie} is preferred over the full eigendecomposition.
\par To our knowledge, our proposed modification to power iteration has not been reported in the literature. Work by Shi and Petrovic to improve power iteration convergence relies on incorporating information about the second eigenvalue, but does not demonstrate the exponential speed of our Algorithm \ref{alg:pie} \cite{modpi}.
\begin{figure}[h]
	\label{fig:first}
	\caption{Histogram of number of iterations needed to reach convergence ($n=100$, 300 matrices). A: real matrices, Algorithm \ref{alg:pi}. B: real matrices, Algorithm \ref{alg:pie}. C: complex matrices, Algorithm \ref{alg:pi}. D: complex matrices, Algorithm \ref{alg:pie}. The distribution of $\log_2$ iterations for convergence in Algorithm \ref{alg:pi} is similar to linear scale for iterations of Algorithm \ref{alg:pie}. Algorithm \ref{alg:pie} reaches convergence within 20 iterations due to the exponential growth of its accuracy. There is long tail to the right in the execution times of both algorithms, which is greatly exaggerated when performing naive power iteration, leading to slow and unpredictable convergence when using Algorithm \ref{alg:pi}.}
	\includegraphics[width=\linewidth]{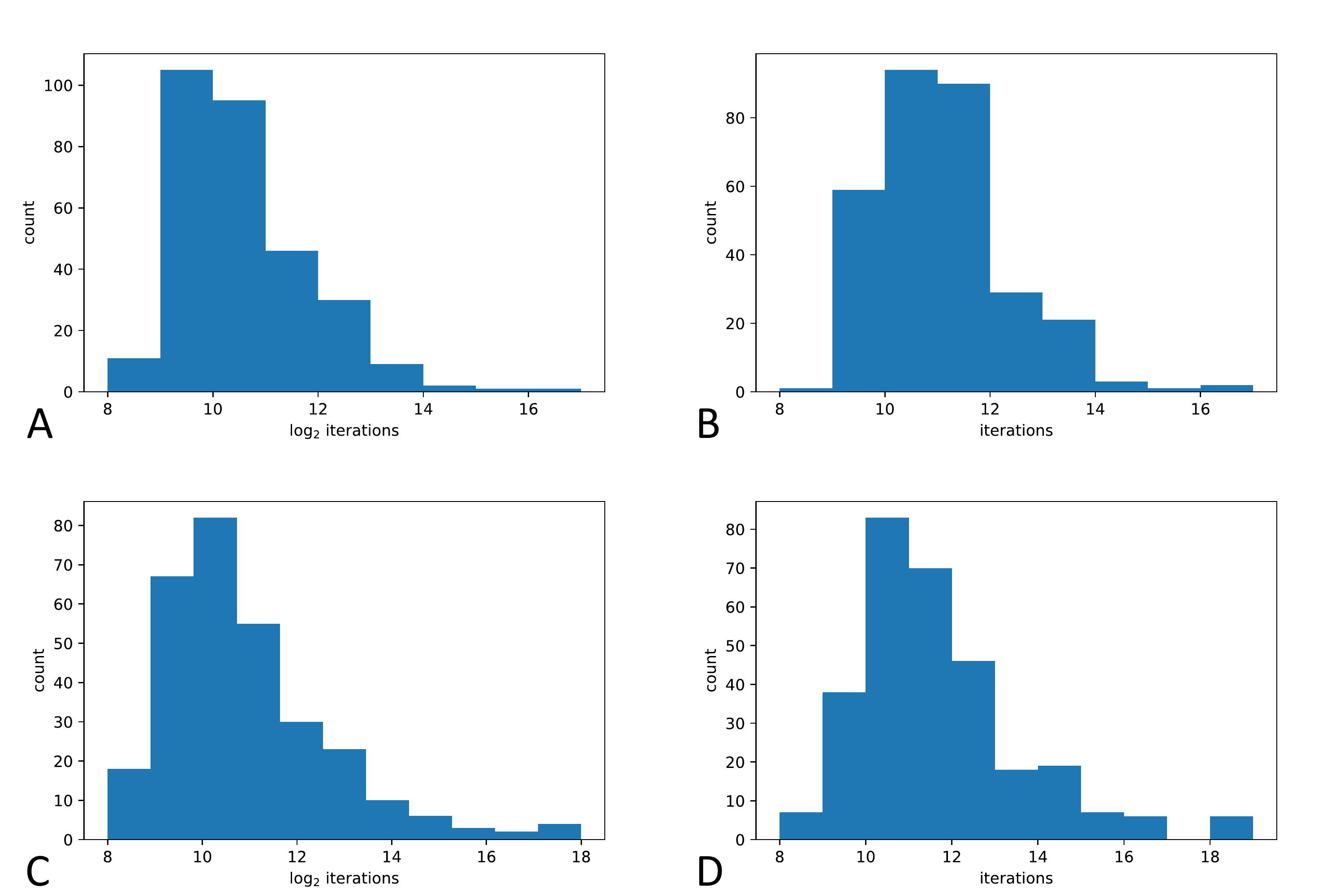}
\end{figure}
\section{Conclusions}
Modifying the power iteration algorithm to take advantage of matrix exponentiation by squaring produces a fast, practical, and robust method of determining the largest eigenvalue and its associated eigenvector for matrices. Our Algorithm \ref{alg:pie} outperforms traditional power iteration in convergence rate and speed on a variety of matrix sizes. We recommend Algorithm \ref{alg:pie} as a drop-in replacement for power iteration due to its favorable performance and simple implementation.
\section{Competing interests}
The authors have no competing interests to disclose.
\section{Funding}
We acknowledge support from the National Institutes for Health (R35 GM134864), National Science Foundation (2040667) and the Passan Foundation. This project is also funded, in part, under a grant with the Pennsylvania Department of Health using Tobacco CURE Funds. The Department specifically disclaims responsibility for any analyses, interpretations or conclusions.
\pagebreak
\bibliographystyle{ieeetr}
\bibliography{bib.bib}
\end{document}